\documentclass[14pt]{article}
\usepackage[left=1in,top=1in,right=1in,bottom=1in,letterpaper]{geometry}

\usepackage{amsmath,amssymb,amsthm}
\usepackage{latexsym,amsfonts,amscd,amsxtra,amstext}

\newcommand{\RR}{\mathbb{R}}
\newcommand{\NN}{\mathbb{N}}
\newcommand{\dom}{{\mathrm{dom}}} 
\newcommand{\prox}{{\mathbf{prox}}}

\DeclareMathOperator*{\Min}{minimize}

\newcommand{\st}{\mbox{subject to}}

\newtheorem{theorem}{Theorem}
\newtheorem{assumption}{Assumption}

\newtheorem{lemma}{Lemma}

\begin{document}

\title{Variants of alternating minimization method with sublinear rates of convergence for convex optimization}

\author{
Hui Zhang\thanks{
College of Science, National University of Defense Technology,
Changsha, Hunan, China, 410073. Email: \texttt{h.zhang1984@163.com}}
\and Lizhi Cheng \thanks{The state key laboratory for high performance computation, and College of Science,
 National University of Defense Technology,
Changsha, Hunan, China, 410073. Email: \texttt{clzcheng@nudt.edu.cn}}
}
\date{\today}

\maketitle

\begin{abstract}
The alternating minimization (AM) method is a fundamental method for minimizing convex functions whose variable consists of two blocks. How to efficiently solve each subproblems when applying the AM method is the most concerned task. In this paper,  we investigate this task and design two new variants of the AM method by borrowing proximal linearized techniques. The first variant is very suitable for the case where half of the subproblems are hard to be solved and the other half can be directly computed. The second variant is designed for parallel computation. Both of them are featured by simplicity at each iteration step. Theoretically, with the help of the proximal operator we first write the new as well as the existing AM variants into uniform expressions, and then prove that they enjoy sublinear rates of convergence under very minimal assumptions.
\end{abstract}

\section{Introduction}
The alternating minimization (AM) method is a fundamental algorithm for solving the following optimization problem:
\begin{equation}\label{generalproblem}
 \Min_{x\in\RR^{n_1},y\in\RR^{n_2}} \Psi(x,y):=f(x)+H(x,y)+g(y),
\end{equation}
where $\Psi(x,y)$ is a convex function.
Starting with a given initial point $(x^0,y^0)$, the AM method generates a sequence $\{(x^k,y^k)\}_{k\in\NN}$ via the scheme
\begin{subequations}
\begin{align}
&x^{k+1}\in\arg\min \{H(x, y^k)+f(x)\}\\
&y^{k+1}\in\arg\min \{H(x^{k+1}, y)+g(y)\}.
\end{align}
\end{subequations}
 In the literature, there exist lots of work concerning its convergence with certain assumptions. To obtain stronger convergence results under more general settings, the recent paper \cite{attouch2010proximal} proposed an augmented alternating minimization (AAM) method by adding proximal terms, that is
 \begin{subequations}\label{AAM}
\begin{align}
&x^{k+1}\in\arg\min \{H(x, y^k)+f(x)+\frac{c_k}{2}\|x-x^k\|^2\}\\
&y^{k+1}\in\arg\min \{H(x^{k+1}, y)+g(y)+\frac{d_k}{2}\|y-y^k\|^2\}
\end{align}
\end{subequations}
where $c_k, d_k$ are positive real numbers. From practical computational perspective, the authors in another recent paper \cite{bolte2013proximal} suggested the proximal alternating linearized minimization (PALM) scheme:
\begin{subequations}\label{PALM}
\begin{align}
&x^{k+1}\in\arg\min_x \{f(x)+\langle x-x^k,\nabla_xH(x^k,y^k)\rangle+\frac{c_k}{2}\|x-x^k\|^2\}\\
&y^{k+1}\in\arg\min_y \{g(y)+\langle y-y^k,\nabla_yH(x^{k+1},y^k)\rangle+\frac{d_k}{2}\|y-y^k\|^2\},
\end{align}
\end{subequations}
By the powerful Kurdyka-{\L}ojasiewicz property, the global convergence of both the augmented AM method and the PALM method were established with very minimal assumptions in \cite{attouch2010proximal,bolte2013proximal}. A remarkable feature of the convergence analysis by the Kurdyka-{\L}ojasiewicz property is to remove the convexity assumption over $\Psi(x, y)$. However, once the convexity is absent, to obtain rates of convergence will become very difficulty or even impossible. Very recently, the author of \cite{beck2014convergence} brought back the convexity and derived the first sublinear $\mathcal{O}(\frac{1}{k})$ rate of convergence for the AM method employed on the general problem \eqref{generalproblem}.

In this paper, we concern both of the computational and theoretical aspects of the AM method. On the computational hand, we follow the proximal linearized technique employed by the PALM method and propose two new variants of the AM method, called AM-variant-I and AM-variant-II. They read as follows respectively:
\begin{subequations}\label{variantI}
\begin{align}
&x^{k+1}\in\arg\min_x \{f(x)+\langle x-x^k,\nabla_xH(x^k,y^k)\rangle+\frac{c_k}{2}\|x-x^k\|^2\}\\
&y^{k+1}\in\arg\min_y \{H(x^{k+1}, y)+g(y)\},
\end{align}
\end{subequations}
and
\begin{subequations}\label{variantII}
\begin{align}
&x^{k+1}\in\arg\min_x \{f(x)+\langle x-x^k,\nabla_xH(x^k,y^k)\rangle+\frac{c_k}{2}\|x-x^k\|^2\}\\
&y^{k+1}\in\arg\min_y \{g(y)+\langle y-y^k,\nabla_yH(x^k,y^k)\rangle+\frac{d_k}{2}\|y-y^k\|^2\},
\end{align}
\end{subequations}
AM-variant-I can be viewed as a hybrid of the original AM method and the PLAM method. The proximal linearized technique is only employed to update the $x$-variable, and the updating of $y$-variable is as same as the original AM method. The idea lying in AM-variant-I is mainly motivated by the iteratively reweighted least square (IRLS) method where the subproblem with respective to (w.r.t.) the $y$-variable can be easily computed but the subproblem w.r.t. the $x$-variable might greatly benefit from proximal linearized techniques. AM-variant-II is very similar to the PALM method. The only difference is that we use $\nabla_yH(x^k,y^k)$ rather than $\nabla_yH(x^{k+1},y^k)$ when update the $y$-variables. The biggest merit of this scheme is that it is very suitable for \textsl{parallel computation}.

On the theoretical hand, motivated by the method in \cite{beck2014convergence}, we summarize a theoretical framework with which we prove that all the AM variants, including the AAM method and the PALM method, have the sublinear $\mathcal{O}(\frac{1}{k})$ rate of convergence under minimal assumptions. The rest of this paper is organized as follows. In section 2, we list some basic properties and formulate all AM-variants into uniform expressions by using the proximal operator. In section 3, we first list all the assumptions that needed for convergence analysis, and then state the main convergence results with a proof sketch. The proof details are postponed to section 6. In section 4, we introduce two applications to show the motivation and advantages of AM-variant-I for some special convex optimization problem. Future work is briefly discussed in section 5.

\section{Mathematical Preliminaries}
In this section, we layout some basic properties about gradient-Lipschitz-continuous functions and the proximal operator, and then formulate all AM-variants into uniform expressions based on the proximal operator.
\subsection{Basic properties}
\begin{lemma}[\cite{nesterov2004introductory}]\label{gradLip}
 Let $h:\RR^n\rightarrow \RR$ be a continuously differentiable function and assume that its gradient $\nabla h$ is Lipschitz continuous with constant $L_h<+\infty$:
 $$\|\nabla h(u)-\nabla h(v)\|\leq L_h\|u-v\|,~~\forall u, v\in\RR^n.$$
 Then,  it holds that
 $$h(u)\leq h(v)+\langle \nabla h(v), u-v\rangle +\frac{L_h}{2}\|u-v\|^2,~~\forall u, v\in \RR^n,$$
 and
 $$\langle \nabla h(u)-h(v), u-v \rangle \leq L_h\|u-v\|^2, ~~\forall u, v\in \RR^n.$$
\end{lemma}

Let $\sigma: \RR^n\rightarrow (-\infty, +\infty]$ be a proper and lower semicontinuous convex function. For given $x\in\RR^n$ and $t>0$, the proximal operator is defined by:
\begin{equation}
\prox_{t}^\sigma(x):=\arg\min_{u\in\RR^n}\{\sigma(u)+\frac{t}{2}\|u-x\|^2\}
\end{equation}

The characterization of the proximal operator given in the following lemma is very important for convergence analysis and will be frequently used later.
\begin{lemma}[\cite{beck2009gradient}]\label{prox}
 Let $\sigma: \RR^n\rightarrow (-\infty, +\infty]$ be a proper and lower semicontinuous convex function. Then
 $$w=\prox_{t}^\sigma(x)$$
 if and only if for any $u\in \dom \sigma$:
 \begin{equation}
 \sigma (u)\geq \sigma(w)+ t\langle x-w, u-w\rangle.
 \end{equation}
\end{lemma}

The next result was established in \cite{bolte2013proximal}; its corresponding result in the convex setting appeared in an earlier paper \cite{beck2009fast}.
\begin{lemma}[\cite{bolte2013proximal}]\label{grad-prox}
Let $h:\RR^n\rightarrow \RR$ be a continuously differentiable function and assume that its gradient $\nabla h$ is Lipschitz continuous with constant $L_h<+\infty$ and let $\sigma ¦Ò:\RR^n\rightarrow \RR$ be a proper and lower semicontinuous function with $\inf_{\RR^n} \sigma>-\infty$. Fix any $t>L_h$. Then, for any $u\in \dom \sigma $ and any $u^+\in \RR^n$ defined by $$u^+\in \prox_t^\sigma(u-\frac{1}{t}\nabla h(u))$$  we have
$$h(u^+)+\sigma(u^+)\leq h(u)+\sigma(u)-\frac{1}{2}(t-L_h)\|u^+-u\|^2.$$
\end{lemma}

\subsection{Uniform AM-variant expressions}
In what follows, we express all AM-variants in a uniform way.
AM-variant-I:
\begin{subequations}\label{pvariantI}
\begin{align}
&x^{k+1}= \prox_{c_k}^f(x^k-\frac{1}{c_k}\nabla_xH(x^k,y^k))\\
&y^{k+1}= \prox_{d_{k+1}}^g(y^{k+1}-\frac{1}{d_{k+1}}\nabla_yH(x^{k+1},y^{k+1})).
\end{align}
\end{subequations}
AM-variant-II:
\begin{subequations}\label{pvariantII}
\begin{align}
&x^{k+1}= \prox_{c_k}^f(x^k-\frac{1}{c_k}\nabla_xH(x^k,y^k))\\
&y^{k+1}= \prox_{d_{k}}^g(y^{k}-\frac{1}{d_{k}}\nabla_yH(x^k,y^k)).
\end{align}
\end{subequations}
The AAM method:
\begin{subequations}\label{paam}
\begin{align}
&x^{k+1}= \prox_{c_k}^f(x^k-\frac{1}{c_k}\nabla_xH(x^{k+1},y^k))\\
&y^{k+1}= \prox_{d_{k}}^g(y^{k}-\frac{1}{d_{k}}\nabla_yH(x^{k+1},y^{k+1})).
\end{align}
\end{subequations}
The PALM method:
\begin{subequations}\label{ppalm}
\begin{align}
&x^{k+1}= \prox_{c_k}^f(x^k-\frac{1}{c_k}\nabla_xH(x^k,y^k))\\
&y^{k+1}= \prox_{d_{k}}^g(y^{k}-\frac{1}{d_{k}}\nabla_yH(x^{k+1},y^k)).
\end{align}
\end{subequations}
These expressions can be easily derived by the first optimal condition and the following basic fact:
$$\prox_{t}^\sigma(\cdot)=[I+\frac{1}{t}\partial \sigma(\cdot)]^{-1}.$$
We omit all the deductions here.

\section{Main results}
For convenience, we let $z=(x,y)\in\RR^{n_1}\times\RR^{n_2}$ and $K(x,y)=f(x)+g(y)$. With these notations, the objective function in problem \eqref{generalproblem} equals to $K(x,y)+H(x,y)$ or $K(z)+H(z)$, and $z^{k}=(x^k,y^k)$.

\subsection{Assumptions and convergence results}
Before stating main results, we make the following basic assumptions throughout the paper:
\begin{assumption}\label{assump1}
The functions $f:\RR^{n_1}\rightarrow (-\infty, +\infty]$ and $g:\RR^{n_2}\rightarrow (-\infty, +\infty]$ are proper and lower semicontinuous convex function satisfying $\inf_{\RR^{n_1}}f>-\infty$ and $\inf_{\RR^{n_2}} g>-\infty$. The function $H(x,y)$ is a continuously differentiable convex function over $\dom f\times \dom g$. The function $\Psi$ satisfies $\inf_{\RR^{n_1}\times\RR^{n_2}}\Psi>-\infty$.
\end{assumption}

\begin{assumption}\label{assump2}
The optimal set of \eqref{generalproblem}, denoted by $Z^*$, is nonempty, and the corresponding optimal value is denoted by $\Psi^*$. The level set $$S=\{z\in \dom f\times \dom g: \Psi(z)\leq \Psi(z^0)\}$$
is compact, where $z^0$ is some given initial point.
\end{assumption}

Besides these two basic assumptions, we need additional assumptions to analyze different AM-variants.

\begin{assumption}\label{assump3}
For any fixed $y$, the gradient $\nabla_x H(x, y)$ is Lipschitz continuous with constant $L_1(y)$:
 $$\|\nabla_x H(x_1, y)-\nabla_x H(x_2, y)\|\leq L_1(y)\|x_1-x_2\|,~~\forall x_1, x_2\in\RR^{n_1}.$$
\end{assumption}

\begin{assumption}\label{assump4}
For any fixed $x$, the gradient $\nabla_y H(x, y)$ is Lipschitz continuous with constant $L_2(x)$:
 $$\|\nabla_y H(x, y_1)-\nabla_y H(x, y_2)\|\leq L_2(x)\|y_1-y_2\|,~~\forall y_1, y_2\in\RR^{n_2}.$$
\end{assumption}

\begin{assumption}\label{assump5}
For any fixed $y$, the gradient $\nabla_y H(x, y)$ w.r.t. the variables $x$ is Lipschitz continuous with constant $L_3(y)$:
 $$\|\nabla_y H(x_1, y)-\nabla_y H(x_2, y)\|\leq L_3(y)\|x_1-x_2\|,~~\forall x_1, x_2\in\RR^{n_1}.$$
\end{assumption}

\begin{assumption}\label{assump6}
For any fixed $x$, the gradient $\nabla_x H(x, y)$ w.r.t. the variables $y$  is Lipschitz continuous with constant $L_4(x)$:
 $$\|\nabla_x H(x, y_1)-\nabla_x H(x, y_2)\|\leq L_4(x)\|y_1-y_2\|,~~\forall y_1, y_2\in\RR^{n_2}.$$
\end{assumption}

\begin{assumption}\label{assump7}
 $\nabla H(z)$ is Lipschitz continuous with constant $L_5$.
\end{assumption}

\begin{assumption}\label{assump8}
For $i=1,2,3, 4$ there exists $\lambda_i^-, \lambda_i^+>0$ such that
\begin{subequations}
\begin{align}
&\inf\{L_1(y^k): k\in\NN\}\geq \lambda_1^-~~~~ \textrm{and}~~~~\inf\{L_2(x^k): k\in\NN\}\geq \lambda_2^- \\
&\sup\{L_1(y^k): k\in\NN\}\leq \lambda_1^+~~~~ \textrm{and}~~~~\sup\{L_2(x^k): k\in\NN\}\leq \lambda_2^+ \\
&\sup\{L_3(y^k): k\in\NN\}\leq \lambda_3^+~~~~ \textrm{and}~~~~\sup\{L_4(x^k): k\in\NN\}\leq \lambda_4^+.
\end{align}
\end{subequations}
\end{assumption}

Define $R=\max_{z\in\RR^{n_1\times n_2}}\max_{z^*\in Z^*}\{\|z-z^*\|: \Psi(z)\leq \Psi(z^0)\};$ then $R\leq +\infty$ from Assumption \ref{assump2}. Now, we are ready to present our main results about the convergence rate of different AM-variants.
\begin{theorem}[The convergence rate of AM-variant-I]\label{thm1}
 Suppose that Assumptions \ref{assump3} and \ref{assump8} hold.
 Take $c_k=\gamma\cdot L_1(y^k)$ with $\gamma>1$ and let $\{(x^k,y^k)\}_{k\in\NN}$ be the sequence generated by AM-variant-I. Then, for all $k\geq 2$
  $$\Psi(z^k)-\Psi^*\leq \max\left\{\left(\frac{1}{2}\right)^{(k-1)/2}(\Psi(z^0)-\Psi^*), \frac{8(\lambda_1^+)^2R^2\gamma^2}{\lambda_1^-(\gamma-1)(k-1)}\right\}.$$
\end{theorem}

\begin{theorem}[The convergence rate of AM-variant-II]
 Suppose that Assumptions \ref{assump3}, \ref{assump4}, \ref{assump7} and \ref{assump8} hold. Denote $\eta=\min\{\frac{\lambda_1^-}{2},\frac{\lambda_2^-}{2}\}$
 Take $c_k=\gamma\cdot L_1(y^k)$ and $d_k=\gamma\cdot L_2(x^k)$ with $\gamma>\frac{L_5}{\eta}$ and let $\{(x^k,y^k)\}_{k\in\NN}$ be the sequence generated by AM-variant-II. Then, for all $k\geq 2$
 $$\Psi(z^k)-\Psi^*\leq \max\left\{\left(\frac{1}{2}\right)^{(k-1)/2}(\Psi(z^0)-\Psi^*), \frac{4\max\{(\lambda_1^+)^2,(\lambda_2^+)^2\}R^2\gamma^2}{(\gamma \eta-L_5)(k-1)}\right\}.$$
\end{theorem}

\begin{theorem}[The convergence rate of the AAM method]
 Suppose that Assumptions \ref{assump6} and \ref{assump8} hold. Let $c_k, d_k$ be positive real numbers such that $\rho_1=\inf\{c_k, d_k: d\in\NN\}>0$ and $\rho_2=\sup\{c_k, d_k: k\in\NN\}<+\infty$. Let $\{(x^k,y^k)\}_{k\in\NN}$ be the sequence generated by the AAM method. Then, for all $k\geq 2$
  $$\Psi(z^k)-\Psi^*\leq \max\left\{\left(\frac{1}{2}\right)^{(k-1)/2}(\Psi(z^0)-\Psi^*), \frac{8R^2(\rho_2+\lambda_4^+)^2}{\rho_1(k-1)}\right\}.$$
\end{theorem}

\begin{theorem}[The convergence rate of the PALM method]
 Suppose that Assumptions \ref{assump3}, \ref{assump4}, \ref{assump5}, and \ref{assump8} hold. Take $c_k=\gamma\cdot L_1(y^k)$ and $d_k=\gamma\cdot L_2(x^{k+1})$ with $\gamma>1$. Let $\{(x^k,y^k)\}_{k\in\NN}$ be the sequence generated by the PALM method. Then, for all $k\geq 2$
   $$\Psi(z^k)-\Psi^*\leq \max\left\{\left(\frac{1}{2}\right)^{(k-1)/2}(\Psi(z^0)-\Psi^*), \frac{8R^2(\lambda_3^++\gamma\max\{\lambda_1^+,\lambda_2^+\})^2}{\min\{\lambda_1^-, \lambda_2^-\}(\gamma-1)(k-1)}\right\}.$$
\end{theorem}

\subsection{Proof sketch}
In the light of \cite{beck2014convergence}, we describe a theoretical framework under which all the theorems stated above can be proved. Assume that a generic algorithm $\mathcal{A}$ generates a sequence $\{z^k\}_{k\in\NN}$ for solving the problem \eqref{generalproblem}. Our aim is to show that
$$\Psi(z^k)-\Psi^*\leq \mathcal{O}(\frac{1}{k}).$$
Our proof mainly consists of two steps:
   \begin{description}
   \item[(a)]  Find a positive constant $\tau_1$ such that
   $$\Psi(z^k)-\Psi(z^{k+1})\geq \tau_1\cdot d(z^k,z^{k+1})^2, ~~k=0, 1, 2, \cdots$$
   where $d(\cdot ,\cdot)$ is some distance function.
   \item[(b)]   Find a positive constant $\tau_2$ such that
   $$\Psi(z^{k+1})-\Psi^* \leq \tau_2\cdot d(z^k,z^{k+1}), ~~k=0, 1, 2, \cdots$$
   \end{description}
Combining these two properties, we can conclude that there exist the positive constant $\alpha=\frac{\tau_1}{\tau_2^2}$ such that
$$(\Psi(z^k)-\Psi^*) -(\Psi(z^{k+1})-\Psi^*) \geq \alpha \cdot (\Psi(z^{k+1})-\Psi^*)^2, ~~k=0, 1, 2, \cdots$$
All the theorems directly follow by invoking the following lemma:
 \begin{lemma}[\cite{beck2014convergence}]
 Let $\{A_k\}_{k\geq 0}$ be a nonnegative sequence of real numbers satisfying
 $$A_k-A_{k+1}\geq \alpha A^2_{k+1}, ~~k=0, 1, 2, \cdots$$
 Then, for any $k\geq 2$,
 $$A_k\leq \max\left\{\left(\frac{1}{2}\right)^{(k-1)/2}A_0, \frac{4}{\alpha(k-1)}\right\}.$$
 \end{lemma}

\section{Applications}
In this part, we first explain our original motivation of proposing AM-variant-I by studying a recent application of the IRLS method; and then we apply  AM-variant-I to solving a composite convex model. We begin with the general problem of minimizing the sum of a continuously differentiable function and sum of norms of affine mappings:
\begin{eqnarray}
\begin{array}{ll}
 & \Min_{x\in\RR^n} \quad s(x)+\sum_{i=1}^m\|A_ix+b_i\| \\
  & \st \quad x\in X,
\end{array}
\end{eqnarray}
where $X$ is a given convex set, $A_i$ and $b_i$ are given matrices and vectors, and $s(x)$ is some continuously differentiable convex function.
This problem was considered and solved in \cite{beck2014convergence} by applying the IRLS method to its smoothed approximation problem:
\begin{eqnarray}
\begin{array}{ll}
 & \Min \quad s(x)+\sum_{i=1}^m\sqrt{\|A_ix+b_i\|^2+\epsilon^2} \\
  & \st \quad x\in X,
\end{array}
\end{eqnarray}
or equivalently, by applying the original AM method to an auxiliary problem
\begin{equation}
 \Min_{x\in\RR^{n},y\in\RR^{m}}  h_\epsilon(x, y)+\delta(x, X)+\delta(y, [\epsilon/2, +\infty)^m),
\end{equation}
where for a given set $Z$ the indicator function $\delta(x, Z)$ is defined by
\begin{eqnarray}
\delta(x, Z)=\left\{\begin{array}{ll}
0 &\textrm{if} ~x\in Z\\
+\infty&\textrm{otherwise}
\end{array} \right.
\end{eqnarray}
and $ h_\epsilon(x, y)=s(x)+\frac{1}{2}\sum_{i=1}^m\left(\frac{\|A_ix+b_i\|^2+\epsilon^2}{y_i} +y_i\right)$. Both of the IRLS and the AM methods generate the same sequence $\{x^k\}_{k\in\NN}$ via the scheme:
\begin{subequations}
\begin{align}
&x^{k+1}\in\arg\min_{x\in X}\{s(x)+ \frac{1}{2}\sum_{i=1}^m\|A_ix+b_i\|^2y_i^k\}\\
&y_i^{k+1}=\frac{1}{\sqrt{\|A_ix^{k+1}+b_i\|^2+\epsilon^2}},~i=1, 2, \cdots, m.
\end{align}
\end{subequations}
In \cite{beck2014convergence}, the author first established the sublinear rate of convergence for the AM method and hence the same convergence result for the IRLS method follows. However, in many case the subproblem of updating $x^k$ is very hard to be solved and even prohibitive for large-scale problems. It is just this drawback  motivating us to propose AM-variant-I. Now, applying AM-variant-I and with some simple calculation, we at once obtain a linearized scheme of the IRLS method, that is
\begin{subequations}\label{linearizedIRLS}
\begin{align}
&x^{k+1}= \mathcal{P}_X\left(x^k-\frac{1}{c_k}(\nabla s(x^k)+ \sum_{i=1}^my_i^kA^T_i(A_ix^k+b_i))\right)\\
&y_i^{k+1}=\frac{1}{\sqrt{\|A_ix^{k+1}+b_i\|^2+\epsilon^2}},~i=1, 2, \cdots, m,
\end{align}
\end{subequations}
where $\mathcal{P}_X$ is the projection operator onto $X$. If $\mathcal{P}_X$ can be easily computed, then the linearized scheme becomes very simple. In addition, its sublinear rate of convergence can be guaranteed by Theorem \ref{thm1}. Nevertheless, we would like to point out that the scheme \eqref{linearizedIRLS} can also be obtained by applying the proximal forward-backward (PFB) method \cite{combettes2005signal} to the following problem
\begin{equation}
\Min_{x\in\RR^n}\quad \delta(x, X)+ s(x)+  \sum_{i=1}^m\sqrt{\|A_ix+b_i\|^2+\epsilon^2}
\end{equation}
and hence can be accelerated into an $\mathcal{O}(\frac{1}{k^2})$-convergent scheme by the Nesterov technique \cite{nesterov2012gradient}. In this sense, AM-variant-I seemly does not bring us more information than the existing methods. Fortunately, the whole fact has not been completely reflected by this motivating example.  To illustrate this, we consider the composite convex model:
\begin{equation}
\Min_{x\in\RR^n}\quad f(x)+g(Ax)
\end{equation}
where $A\in \RR^{m\times n}$ and the proximal operators of $f$ and $g$ can be easily computed. In \cite{beck2014convergence}, the author applied the AM method to its auxiliary problem
\begin{equation}
\Min_{x\in\RR^{n},y\in\RR^{m}}\quad f(x)+g(y)+\frac{\rho}{2}\|Ax-y\|^2
\end{equation}
and obtain the following scheme:
\begin{subequations}
\begin{align}
&x^{k+1}\in\arg\min\{f(x)+ \frac{\rho}{2}\|Ax-y^k\|^2\}\\
&y^{k+1}\in\arg\min\{g(y)+ \frac{\rho}{2}\|Ax^{k+1}-y\|^2\}.
\end{align}
\end{subequations}
Because the entries of vector $x$ are coupled by $Ax$, the updating of $x^k$ is usually very hard for large-scale problems. AM-variant-I fits into this problem and generates the following simple scheme:
\begin{subequations}
\begin{align}
&x^{k+1}=\prox_{c_k}^f(x^k-\frac{\rho}{c_k}A^T(Ax^k-y^k))\\
&y^{k+1}=\prox_{\rho}^g(Ax^{k+1}),
\end{align}
\end{subequations}
where $\rho>0$ is the penalty parameter and $c_k=\gamma\rho\|AA^T\|$ with $\gamma>1$ is the step size parameter. Its sublinear rate of convergence is then guaranteed by Theorem \ref{thm1}.

\section{Discussion}
In this paper, we discuss a group of variants of the AM method and derive sublinear rates of convergence under very minimal assumptions.
Although we restrict our attention onto convex optimization problems, these variants for nonconvex optimization problems might obtain computational advantages over the AM method as well. Because our theory is limited to convex optimization, the convergence of AM-variant-I and AM-variant-II for general cases is unclear at present. In future pursuit, we will analyze the convergence of AM-variant-I and AM-variant-II under nonconvex setting.


\section{Proof details}

\subsection{Proof of Theorem 1}
\textbf{Step 1: prove the property (a).} Denote $z^{k+\frac{1}{2}}=(x^{k+1},y^k)$. On one hand, by invoking Lemma \ref{grad-prox} we derive that
\begin{subequations}
\begin{align}
 \Psi(z^k)-\Psi(z^{k+\frac{1}{2}})&= f(x^k)+H(x^k, y^k)-(f(x^{k+1})+  H(x^{k+1}, y^k))\nonumber\\
  &\geq  \frac{c_k-L_1(y^k)}{2} \|x^{k+1}-x^k\|^2=\frac{\gamma-1}{2}L_1(y^k) \|x^{k+1}-x^k\|^2\label{iter1r}\\
  &\geq  \frac{(\gamma-1)\lambda_1^-}{2} \|x^{k+1}-x^k\|^2\label{ineq3}.
\end{align}
\end{subequations}
On the other hand, since $y^{k+1}$ minimizes the objective $H(x^{k+1},y) +g(y)$, it holds that
\begin{equation} \label{ineq4}
 \Psi(z^{k+\frac{1}{2}})-\Psi(z^{k+1}) = H(x^{k+1},y^k) +g(y^k) -(H(x^{k+1},y^{k+1}) +g(y^{k+1}))\geq 0.
\end{equation}
 By summing the above two inequalities, we obtain
 \begin{equation}
 \Psi(z^k)-\Psi(z^{k+1})\geq \frac{(\gamma-1)\lambda_1^-}{2} \|x^{k+1}-x^k\|^2, ~~k=0, 1, 2, \cdots
\end{equation}

\textbf{Step 2: prove the property (b).} By Lemma \ref{gradLip}, we have that
\begin{align}
 H(z^{k+\frac{1}{2}})-H(z^*)&= H(x^{k+1},y^k)-H(z^*)\nonumber\\
  &\leq  H(x^{k},y^k)+\langle \nabla_x H(x^k,y^k), x^{k+1}-x^k\rangle +\frac{L_1(y^k)}{2}\|x^{k+1}-x^k\|^2-H(z^*)\nonumber\\
  & =H(z^k) + \langle \nabla  H(z^k), z^{k+\frac{1}{2}}-z^k \rangle +\frac{c_k}{2\gamma} \|x^{k+1}-x^k\|^2 -H(z^*).
\end{align}
By the convexity of $H(z)$, it follows that $H(z^k)-H(z^*)\leq \langle \nabla H(z^k), z^k-z^*\rangle$. Thus,
 \begin{equation}\label{ineq1}
 H(z^{k+\frac{1}{2}})-H(z^*)\leq \langle \nabla  H(z^k), z^{k+\frac{1}{2}}-z^* \rangle +\frac{c_k}{2\gamma} \|x^{k+1}-x^k\|^2.
\end{equation}
Recall that
$x^{k+1}= \prox_{c_k}^f(x^k-\frac{1}{c_k}\nabla_xH(x^k,y^k))$ and
$y^{k}= \prox_{d_{k}}^f(y^{k}-\frac{1}{d_{k}}\nabla_yH(x^{k},y^{k}))$. Applying Lemma \ref{prox} to them, we obtain that
  \begin{equation}
 f(x^*)\geq f(x^{k+1})+c_k\langle x^k-x^{k+1}, x^*-x^{k+1}\rangle +\langle \nabla_x H(x^k,y^k), x^{k+1}-x^*\rangle
\end{equation}
and
 \begin{equation}
g(y^*)\geq g(y^k)+\langle \nabla_y H(x^k,y^k), y^{k}-y^*\rangle.
\end{equation}
By summing the above two inequalities, we obtain
\begin{equation}\label{ineq2}
K(z^*)\geq K(z^{k+\frac{1}{2}})+ c_k\langle x^k-x^{k+1}, x^*-x^{k+1}\rangle + \langle \nabla H(z^k), z^{k+\frac{1}{2}}-z^*\rangle.
\end{equation}
Combining inequalities \eqref{ineq1} and \eqref{ineq2} and noticing $\gamma>1$, we have that
\begin{align}
 \Psi(z^{k+\frac{1}{2}})-\Psi^*&\leq  \frac{c_k}{2\gamma} \|x^{k+1}-x^k\|^2 - c_k\langle x^k-x^{k+1}, x^*-x^{k+1}\rangle \nonumber\\
  &\leq  c_k\|x^{k+1}-x^k\|^2 - c_k\langle x^k-x^{k+1}, x^*-x^{k+1}\rangle \nonumber\\
  & \leq c_k\langle x^{k+1}-x^k, x^*-x^k\rangle\nonumber\\
  &\leq \gamma \lambda_1^+\| x^{k+1}-x^k\|\cdot\|x^*-x^k\|.
\end{align}
From inequalities \eqref{ineq3} and \eqref{ineq4}, it follows that
 \begin{equation}\label{ineq5}
\Psi(z^{k})-\Psi^*\geq \Psi(z^{k+\frac{1}{2}})-\Psi^* \geq \Psi(z^{k+1})-\Psi^*\geq 0.
\end{equation}
Then, $\|x^*-x^k\|\leq \|z^*-z^k\|\leq R$ follows from \eqref{ineq5} and the definition of the constant $R$. Finally, we get
 \begin{equation}
\Psi(z^{k+1})-\Psi^* \leq  R \gamma \lambda_1^+\| x^{k+1}-x^k\|.
\end{equation}

Denote $\tau_1=\frac{(\gamma-1)\lambda_1^-}{2}$ and $\tau_2=R \gamma \lambda_1^+$; then $\alpha=\frac{(\gamma-1)\lambda_1^-}{2R^2\gamma^2(\lambda_1^+)^2}$. This completes the proof of Theorem 1.

\subsection{Proof of Theorem 2}
\textbf{Step 1: prove the property (a).} Since $x^{k+1}$ and $y^{k+1}$ are the minimizers to the subproblems in \eqref{variantII} respectively, we get that
  \begin{equation}
 f(x^k)\geq f(x^{k+1})+ \frac{c_k}{2}\|x^{k+1}-x^k\|^2+ \langle x^{k+1}- x^k, \nabla_x H(x^k,y^k) \rangle
\end{equation}
and
 \begin{equation}
 g(x^k)\geq g(x^{k+1})+ \frac{d_k}{2}\|y^{k+1}-y^k\|^2+ \langle y^{k+1}- y^k, \nabla_y H(x^k,y^k) \rangle.
\end{equation}
By summing the above two inequalities, we obtain
 \begin{equation}
 K(z^k)\geq K(z^{k+1})+ \frac{c_k}{2}\|x^{k+1}-x^k\|^2+ \frac{d_k}{2}\|y^{k+1}-y^k\|^2+ \langle z^{k+1}- z^k, \nabla H(z^k) \rangle.
\end{equation}
By the convexity of $H(z)$, it follows that
\begin{equation}H(z^k)-H(z^{k+1})\geq \langle \nabla H(z^{k+1}), z^k-z^{k+1}\rangle.\end{equation}
By summing the above two inequalities, we get
 \begin{equation}\label{ineq6}
\Psi(z^k)-\Psi(z^{k+1})\geq \frac{c_k}{2}\|x^{k+1}-x^k\|^2+ \frac{d_k}{2}\|y^{k+1}-y^k\|^2 + \langle z^k-z^{k+1}, \nabla H(z^{k+1})-\nabla H(z^k) \rangle.
\end{equation}
By Assumption \ref{assump7}, it holds that
\begin{equation}\label{ineq7}
\langle z^{k+1}-z^k, \nabla H(z^{k+1})-\nabla H(z^k)\rangle\leq L_5\|z^{k+1}-z^k\|^2.
\end{equation}
Since $c_k=\gamma\cdot L_1(y^k)\geq \gamma \lambda_1^-$ and $d_k=\gamma\cdot L_2(x^k)\geq \gamma \lambda_2^-$, from the expression of $\eta$ we have that
\begin{equation}\label{ineq8}
\frac{c_k}{2}\|x^{k+1}-x^k\|^2+ \frac{d_k}{2}\|y^{k+1}-y^k\|^2\geq \gamma\eta \|z^{k+1}-z^k\|^2.
\end{equation}
Now, combining inequalities \eqref{ineq6}, \eqref{ineq7}, and \eqref{ineq8}, we finally get that
\begin{equation}
\Psi(z^k)-\Psi(z^{k+1})\geq  (\gamma\eta-L_5) \|z^{k+1}-z^k\|^2.
\end{equation}
Let $\tau_1=\gamma\eta-L_5$; then it must be positive since $\gamma>\frac{L_5}{\eta}$.

\textbf{Step 2: prove the property (b).} By Assumption \ref{assump7} and Lemma \ref{gradLip}, we have
\begin{equation}
H(z^{k+1})-H(z^*)\leq H(z^k) +\langle \nabla H(z^k), z^{k+1}-z^k\rangle +\frac{L_5}{2}\|z^{k+1}-z^k\|^2-H(z^*).
\end{equation}
The convexity of $H(z)$ implies $H(z^k) -H(z^*)\leq \langle \nabla H(z^k), z^k-z^* \rangle$. Thus, we get
\begin{equation}\label{ineq9}
H(z^{k+1})-H(z^*)\leq \langle \nabla H(z^k), z^{k+1}-z^*\rangle +\frac{L_5}{2}\|z^{k+1}-z^k\|^2.
\end{equation}
Applying Lemma \ref{prox} to \eqref{pvariantII}, we obtain that
  \begin{equation}
 f(x^*)\geq f(x^{k+1})+c_k\langle x^k-x^{k+1}, x^*-x^{k+1}\rangle +\langle \nabla_x H(x^k,y^k), x^{k+1}-x^*\rangle
\end{equation}
and
 \begin{equation}
g(y^*)\geq g(y^{k+1})+d_k\langle y^k-y^{k+1}, y^*-y^{k+1}\rangle + \langle \nabla_y H(x^k,y^k), y^{k+1}-y^*\rangle.
\end{equation}
By summing the above two inequalities, we get
\begin{equation} \label{ineq10}
K(z^*)\geq K(z^{k+1})+ c_k\langle x^k-x^{k+1}, x^*-x^{k+1}\rangle +d_k\langle y^k-y^{k+1}, y^*-y^{k+1}\rangle  + \langle \nabla H(z^k), z^{k+1}-z^*\rangle.
\end{equation}
Combining inequalities \eqref{ineq9} and \eqref{ineq10}, we have
\begin{equation}
\Psi(z^{k+1})-\Psi^*\leq \frac{L_5}{2}\|z^{k+1}-z^k\|^2 - c_k\langle x^k-x^{k+1}, x^*-x^{k+1}\rangle - d_k\langle y^k-y^{k+1}. y^*-y^{k+1}\rangle.
\end{equation}
By the setting of $\gamma>\frac{L_5}{\eta}$, we can deduce that $c_k\geq 2L_5$ and $ d_k\geq 2L_5$. Thus,
\begin{align}
\Psi(z^{k+1})-\Psi^* \leq   &c_k(\|x^{k+1}-x^k\|^2-\langle x^k-x^{k+1}, x^*-x^{k+1}\rangle) \nonumber\\
&+d_k(\|y^{k+1}-y^k\|^2-\langle y^k-y^{k+1}, y^*-y^{k+1}\rangle)\nonumber\\
=& c_k\langle x^{k+1}-x^k, x^*-x^k\rangle +d_k\langle y^{k+1}-y^k, y^*-y^k\rangle \nonumber\\
\leq &c_k\|x^{k+1}-x^k\|\cdot\|x^*-x^k\| +d_k\|y^{k+1}-y^k\|\cdot\|y^*-y^k\|\nonumber\\
\leq & \gamma \lambda_1^+ \|x^{k+1}-x^k\|\cdot\|x^*-x^k\| +\gamma\lambda_2^+\|y^{k+1}-y^k\|\cdot\|y^*-y^k\|\nonumber\\
\leq &\gamma\max\{\lambda_1^+,\lambda_2^+\} \|z^{k+1}-z^k\|\cdot\|z^*-z^k\|\label{rela4},
\end{align}
where the last relationship follows from the Cauchy-Schwartz inequality. In \textbf{step 1}, we have shown that $\Psi(z^k)\geq \Psi(z^{k+1})$ and hence $\|z^*-z^k\|\leq R$ from the definition of the constant $R$. Thus, we finally get
\begin{equation}
\Psi(z^{k+1})-\Psi^* \leq \gamma\max\{\lambda_1^+,\lambda_2^+\} R\|z^{k+1}-z^k\|.
\end{equation}
Let $\tau_2=\gamma\max\{\lambda_1^+,\lambda_2^+\}R$; then the positive $\alpha=\tau_1/\tau_2^2$ exists. This completes the proof.

\subsection{Proof of Theorem 3}
\textbf{Step 1: prove the property (a).} Since $x^{k+1}$ and $y^{k+1}$ are the minimizers to the subproblems in \eqref{paam} respectively, we get that
  \begin{equation}
 H(x^k,y^k)+f(x^k)\geq  H(x^{k+1}, y^k)+ f(x^{k+1})+ \frac{c_k}{2}\|x^{k+1}-x^k\|^2
\end{equation}
and
 \begin{equation}
 H(x^{k+1}, y^k)+g(y^k)\geq  H(x^{k+1}, y^{k+1})+g(y^{k+1}) + \frac{d_k}{2}\|y^{k+1}-y^k\|^2.
\end{equation}
By summing the above two inequalities, we obtain
 \begin{equation}
\Psi(z^k)\geq \Psi(z^{k+1})+ \frac{c_k}{2}\|x^{k+1}-x^k\|^2+ \frac{d_k}{2}\|y^{k+1}-y^k\|^2.
\end{equation}
Since $\rho_1 =\inf\{c_k, d_k:  k\in\NN\}>0$, we have
\begin{equation}
\frac{c_k}{2}\|x^{k+1}-x^k\|^2+ \frac{d_k}{2}\|y^{k+1}-y^k\|^2\geq \frac{\rho_1}{2} \|z^{k+1}-z^k\|^2.
\end{equation}
Thus, we finally get
\begin{equation}
\Psi(z^k)-\Psi(z^{k+1})\geq  \frac{\rho_1}{2} \|z^{k+1}-z^k\|^2.
\end{equation}

\textbf{Step 2: prove the property (b).}
Applying Lemma \ref{prox} to \eqref{paam}, we obtain that
  \begin{equation}
 f(x^*)\geq f(x^{k+1})+c_k\langle x^k-x^{k+1}, x^*-x^{k+1}\rangle +\langle \nabla_x H(x^{k+1},y^k), x^{k+1}-x^*\rangle
\end{equation}
and
 \begin{equation}
g(y^*)\geq g(y^{k+1})+d_k\langle y^k-y^{k+1}, y^*-y^{k+1}\rangle + \langle \nabla_y H(x^{k+1},y^{k+1}), y^{k+1}-y^*\rangle.
\end{equation}
By summing the above two inequalities and letting $$\tilde{\nabla}_{1,k}=( \nabla_x H(x^{k+1},y^k), \nabla_y H(x^{k+1},y^{k+1})),$$ we obtain
\begin{equation}\label{relaapp1}
K(z^*)\geq K(z^{k+1})+ c_k\langle x^k-x^{k+1}, x^*-x^{k+1}\rangle +d_k\langle y^k-y^{k+1}, y^*-y^{k+1}\rangle  + \langle \tilde{\nabla}_{1,k}, z^{k+1}-z^*\rangle.
\end{equation}
By the convexity of $H(z)$, it follows that
\begin{equation}\label{convapp1}
H(z^{k+1})-H(z^*)\leq \langle \nabla H(z^{k+1}), z^{k+1}-z^*\rangle.
\end{equation}
Thus, combining \eqref{relaapp1} and \eqref{convapp1} yields
\begin{align*}
\Psi(z^{k+1})-\Psi^*\leq  &\langle \nabla_x H(x^{k+1},y^{k+1})
-\nabla_x H(x^{k+1},y^k), x^{k+1}-x^*\rangle \\
&-c_k\langle x^k-x^{k+1}, x^*-x^{k+1}\rangle - d_k\langle y^k-y^{k+1}, y^*-y^{k+1}\rangle.
\end{align*}
By Assumptions \ref{assump6} and \ref{assump8}, we deduce that
\begin{subequations}
\begin{align}
&\langle\nabla_x H(x^{k+1},y^{k+1})-\nabla_x H(x^{k+1},y^k), x^{k+1}-x^*\rangle\nonumber \\
 \leq &\|\nabla_x H(x^{k+1},y^{k+1})-\nabla_x H(x^{k+1},y^k)\|\cdot \|x^{k+1}-x^*\|,\\
\leq &L_4(x^{k+1})\|y^{k+1}-y^{k}\|\cdot\|x^{k+1}-x^*\|\\
\leq &\lambda_4^+\|y^{k+1}-y^{k}\|\cdot\|x^{k+1}-x^*\|.
\end{align}
\end{subequations}
By the Cauchy-Schwartz inequality and the notation $\rho_2=\sup\{c_k,d_k: k\in\NN\}$, we have
\begin{equation}
-c_k\langle x^k-x^{k+1}, x^*-x^{k+1}\rangle - d_k\langle y^k-y^{k+1}, y^*-y^{k+1}\rangle\leq \rho_2\|z^k-z^{k+1}\|\cdot\| z^*-z^{k+1}\|.
\end{equation}
As same as before, it holds that $\| z^*-z^{k+1}\|\leq R$. Hence
\begin{subequations}
\begin{align}
\Psi(z^{k+1})-\Psi^*\leq &\lambda_4^+\|y^{k+1}-y^{k}\|\cdot\|x^{k+1}-x^*\|+ \rho_2\|z^k-z^{k+1}\|\cdot\| z^*-z^{k+1}\|\\
\leq &(\lambda_4^++\rho_2)\|z^k-z^{k+1}\|\cdot\| z^*-z^{k+1}\|\leq R(\lambda_4^++\rho_2)\|z^k-z^{k+1}\|.
\end{align}
\end{subequations}
This completes the proof.

\subsection{Proof of Theorem 4}
\textbf{Step 1: prove the property (a).} The following proof appeared in \cite{bolte2013proximal}. For completion, we include it here. Applying Lemma \ref{grad-prox} to the scheme \eqref{ppalm}, we derive that
\begin{subequations}
\begin{align}
 H(x^{k+1}, y^k) +f(x^{k+1})&\leq H(x^k, y^k) +f(x^k)-\frac{1}{2}(c_k-L_1(y^k))\|x^{k+1}-x^k\|^2 \\
  &= H(x^k, y^k) +f(x^k)-\frac{\gamma-1}{2}L_1(y^k)\|x^{k+1}-x^k\|^2\label{iter1r}
\end{align}
\end{subequations}
and
\begin{subequations}
\begin{align}
 H(x^{k+1}, y^{k+1}) +g(x^{k+1})&\leq H(x^{k+1}, y^k) +g(y^k)-\frac{1}{2}(d_k-L_2(x^k))\|y^{k+1}-y^k\|^2 \\
  &= H(x^{k+1}, y^k) +g(y^k)-\frac{\gamma-1}{2}L_2(x^k)\|y^{k+1}-y^k\|^2\label{iter1r}.
\end{align}
\end{subequations}
By summing the above two inequalities, we obtain that
 \begin{equation}
\Psi(z^k)\geq \Psi(z^{k+1})+ \frac{\gamma-1}{2}L_1(y^k)\|x^{k+1}-x^k\|^2+ \frac{\gamma-1}{2}L_2(x^k)\|y^{k+1}-y^k\|^2.
\end{equation}
By Assumption \ref{assump8}, we finally get that
 \begin{equation}
\Psi(z^k)- \Psi(z^{k+1})\geq \frac{\gamma-1}{2}\min\{\lambda_1^-, \lambda_2^-\}\|z^{k+1}-z^k\|^2.
\end{equation}

\textbf{Step 2: prove the property (b).}
Applying Lemma \ref{prox} to \eqref{ppalm}, we obtain that
  \begin{equation}
 f(x^*)\geq f(x^{k+1})+c_k\langle x^k-x^{k+1}, x^*-x^{k+1}\rangle +\langle \nabla_x H(x^{k},y^k), x^{k+1}-x^*\rangle
\end{equation}
and
 \begin{equation}
g(y^*)\geq g(y^{k+1})+d_k\langle y^k-y^{k+1}, y^*-y^{k+1}\rangle + \langle \nabla_y H(x^{k+1},y^{k}), y^{k+1}-y^*\rangle.
\end{equation}
By summing the above two inequalities and letting $$\tilde{\nabla}_{2,k}=( \nabla_x H(x^{k},y^k), \nabla_y H(x^{k+1},y^k)),$$ we obtain that
\begin{equation}\label{rela3}
K(z^*)\geq K(z^{k+1})+ c_k\langle x^k-x^{k+1}, x^*-x^{k+1}\rangle +d_k\langle y^k-y^{k+1}, y^*-y^{k+1}\rangle  + \langle \tilde{\nabla}_{2,k}, z^{k+1}-z^*\rangle.
\end{equation}
 By Assumption \ref{assump4} and Lemma \ref{gradLip}, we have that
 \begin{subequations}\label{rela1}
\begin{align}
 H(z^{k+1})-H(z^{k+\frac{1}{2}})& = H(x^{k+1}, y^{k+1})- H(x^{k+1}, y^{k})\\
  &\leq   \langle \nabla_y H(x^{k+1},y^{k}), y^{k+1}-y^k\rangle+  \frac{d_k}{2\gamma}\|y^{k+1}-y^k\|^2\label{iter1r}.
\end{align}
\end{subequations}
By Assumption \ref{assump3} and Lemma \ref{gradLip} and the convexity of $H(z)$, we derive that
 \begin{subequations}\label{rela2}
\begin{align}
 H(z^{k+\frac{1}{2}})-H(z^*) & \leq  H(x^{k},y^k)+\langle \nabla_x H(x^k,y^k), x^{k+1}-x^k\rangle +\frac{L_1(y^k)}{2}\|x^{k+1}-x^k\|^2-H(z^*) \\
  & \leq \langle \nabla  H(z^k), z^k-z^* \rangle +\langle \nabla_x H(x^k,y^k), x^{k+1}-x^k\rangle  +\frac{c_k}{2\gamma} \|x^{k+1}-x^k\|^2.
\end{align}
\end{subequations}
By summing inequalities \eqref{rela1} and \eqref{rela2}, we get
\begin{equation}
H(z^{k+1})-H(z^*) \leq \langle \nabla  H(z^k), z^k-z^* \rangle +  \langle \tilde{\nabla}_{2,k}, z^{k+1}-z^k\rangle+ \frac{c_k}{2\gamma} \|x^{k+1}-x^k\|^2+\frac{d_k}{2\gamma}\|y^{k+1}-y^k\|^2.
\end{equation}
Together with \eqref{rela3} and utilizing Assumption \ref{assump5} and the fact $\gamma>1$, we derive that
\begin{subequations}
\begin{align}
\Psi(z^{k+1})-\Psi^*  \leq  &\langle \nabla  H(z^k), z^k-z^* \rangle +  \langle \tilde{\nabla}_{2,k}, z^{k+1}-z^k\rangle-\langle \tilde{\nabla}_{2,k}, z^{k+1}-z^*\rangle+\frac{c_k}{2\gamma} \|x^{k+1}-x^k\|^2 \nonumber\\
  &+\frac{d_k}{2\gamma}\|y^{k+1}-y^k\|^2-c_k\langle x^k-x^{k+1}, x^*-x^{k+1}\rangle -d_k\langle y^k-y^{k+1}, y^*-y^{k+1}\rangle\\
 = &\langle \nabla_y H(x^k,y^k)-\nabla_y H(x^{k+1},y^{k}), y^*-y^k \rangle +\frac{c_k}{2\gamma} \|x^{k+1}-x^k\|^2 \nonumber\\
  &+\frac{d_k}{2\gamma}\|y^{k+1}-y^k\|^2-c_k\langle x^k-x^{k+1}, x^*-x^{k+1}\rangle -d_k\langle y^k-y^{k+1}, y^*-y^{k+1}\rangle\\
\leq  & L_3(y^k)\|x^{k+1}-x^k \|\cdot\|y^*-y^k\| +c_k \|x^{k+1}-x^k\|^2 +d_k\|y^{k+1}-y^k\|^2\nonumber\\
  &-c_k\langle x^k-x^{k+1}, x^*-x^{k+1}\rangle -d_k\langle y^k-y^{k+1}, y^*-y^{k+1}\rangle \nonumber\\
 \leq & \lambda_3^+\|z^{k+1}-z^k \|\cdot\|z^*-z^k\|+\gamma\max\{\lambda_1^+,\lambda_2^+\} \|z^{k+1}-z^k\|\cdot\|z^*-z^k\|  \label{rela5}\\
 \leq & (\lambda_3^++\gamma\max\{\lambda_1^+,\lambda_2^+\})\|z^{k+1}-z^k \|\cdot\|z^*-z^k\|\\
 \leq & R(\lambda_3^++\gamma\max\{\lambda_1^+,\lambda_2^+\})\|z^{k+1}-z^k \|,
\end{align}
\end{subequations}
where the second term in \eqref{rela5} can be derived as that in \eqref{rela4}, and the last relationship follows from $\| z^*-z^{k+1}\|\leq R$. This completes the proof.

\section*{Acknowledgements}
The work of L. Cheng is supported by NSF Grants No.61271014 and No.61072118.


\end{document}